\documentclass[smallextended]{svjour3}

\usepackage[paperwidth=6.5in, paperheight=9.5in]{geometry}
\usepackage[latin1]{inputenc}
\usepackage{mathptmx}
\usepackage{amssymb}
\usepackage[firstpage]{draftwatermark}

\newtheorem{te}{Theorem}[section]
\newtheorem{la}{Lemma}[section]
\newtheorem{co}{Corolary}[section]
\newtheorem{po}{Proposition}[section]
\newtheorem{de}{Definition}[section]
\newtheorem{ex}{Example}[section]
\newtheorem{re}{Remark}[section]
\newtheorem{no}{Notation}[section]

\newcommand{\bpr}{\begin{proof}}
\newcommand{\epr}{\end{proof}}
\newcommand{\bte}{\begin{te}}
\newcommand{\ete}{\end{te}}
\newcommand{\bla}{\begin{la}}
\newcommand{\ela}{\end{la}}
\newcommand{\bco}{\begin{co}}
\newcommand{\eco}{\end{co}}
\newcommand{\bpo}{\begin{po}}
\newcommand{\epo}{\end{po}}
\newcommand{\bde}{\begin{de}}
\newcommand{\ede}{\end{de}}
\newcommand{\bex}{\begin{ex}}
\newcommand{\eex}{\end{ex}}
\newcommand{\bre}{\begin{re}}
\newcommand{\ere}{\end{re}}
\newcommand{\bno}{\begin{no}}
\newcommand{\eno}{\end{no}}

\begin{document}

\title{A quantitative approach to semantic informativity}

\author{Anderson de Ara\'ujo}

\institute{Anderson de Ara\'ujo \at
              Federal University of ABC (UFABC) \\
							Center for Natural and Human Sciences (CCNH)\\
							S\~ao Bernardo do Campo, SP, Brazil \\
              \email{anderson.araujo@ufabc.edu.br} 
}

\date{}

\maketitle

\begin{abstract}
This article shows a form of measuring semantic informativity of deductions. Dynamic concepts of complexity and relevancy are presented according to explicit definitions of insertion and deletion on databases. Hence, with respect to finite databases, it solves Bar-Hillel-Carnap paradox and Hintikka' scandal of deduction.
\end{abstract}

\section{Introduction}

It seems plausible to believe we can get new information by reasoning deductively. Concerning this subject, Aristotle said that ``every belief comes either through syllogism or from induction'' (\cite{aristotle1989}, Book 2, Part 23). There is, however, a problem with that basic belief about informativity of deductions:

\begin{enumerate}
\item A deduction $\phi_{0},\phi_{1}\ldots,\phi_{n}$ is valid if, and only
if, the conjunction of its premisses, says $\phi_{0},\ldots,\phi_{k}$,
implies its conclusion, $\phi_{n}$. 
\item In this case, its \emph{associated conditional} $(\phi_{0}\wedge\cdots\wedge\phi_{k})\to\phi_{n}$
is a tautology.
\item Tautologies are propositions without information.
\item Therefore, deductions are uninformative.
\end{enumerate}

Hintikka called this problem \emph{scandal of deduction} \cite{hintikka1970}. It is a scandal, not merely because it contradicts Aristotle's maxim that deductions are important for obtaining beliefs, but, mainly, in virtue of the fact that we actually get information via deductions. In what sense, then, does deductive reasoning give us new information?

For Bar-Hillel-Carnap's classical theory of semantic information \cite{bar-hillel1963}, the information associated to deductions is semantic, it is an information about the meaning of propositions. To be precise, Bar-Hillel-Carnap's theory was designed to measure the amount of semantic information of propositions, not deductions, but, considering the associated conditionals of deductions, it can be extended to deductions too. Bar-Hillel-Carnap's theory is a complement of Shannon's theory of communication \cite{shannon1948}, which treats syntactic information, and so it relies on probability theory - indeed on Carnap's inductive logic. For this reason, it has a very counter-intuitive consequence, namely, contradictions have infinite semantic information, whereas tautologies have null information. Floridi refers this undesirable consequence as the \emph{Bar-Hillel-Carnap paradox} (\cite{floridi2004}). Because tautologies have null information in Bar-Hillel-Carnap's theory, the scandal of deduction has no solution in it.

Hintikka \cite{hintikka1970} displayed a theory of surface information and depth information of propositions, to solve the scandal of deduction in polyadic first-order languages. Analyzing the associated conditional of deductions, Hintikka's theory provided an undecidable measure of semantic informativity that is unpractical (Cf. \cite{sequoiah2008} for details). To overcome that, in \cite{dagostino2009} D'Agostino and Floridi tried to establish a practical distinction between informative and uninformative deductions, using notions from computational complexity. D'Agostino-Floridi's work relies, in turn, on Floridi's \emph{theory of strong semantic information} (Cf. \cite{floridi2011}). According to Floridi' strong theory, semantic information is true well-formed data. From this standpoint and using situation logic, he was able to explain in \cite{floridi2004} why some propositions are more informative than others, and so a solution to Bar-Hillel-Carnap paradox was given, but not to Hintikka's scandal of deduction. This the reason for D'Agostino-Floridi latter work \cite{dagostino2009}, which focused on propositional logic. Nonetheless, D'Agostino-Floridi's approach is, in principle, limited, because their approach can be applied, at best, to monadic first-order languages, in virtue of the restriction of computational complexity to decidable sets.

The aim of this paper is to extend Floridi's original strong perspective to general first-order deductions, without assuming truth as part of the \emph{definiens} of semantic information. This may seem paradoxical, but we will show that it is not. Let us sketch an explanation for that.

We will study deductions into a direct semantic way. This means two things. First, deductions will be analyzed inmselves, and not by their associated conditionals. Second, deductions will be studied in terms of the structures that interpret them, which are called \emph{databases} (Cf. \cite{kroenke2007}). As databases are finite structures and there is, in general, no complete deductive first-order logical system for finite structures (Cf. \cite{ebbinghaus1999}), our semantic approach to deductions is technically well-justified. Indeed, due to the incompleteness of first-order logic for databases, syntactic approaches will be necessarily incomplete. This is our methodological strategy.

In reference to conceptual strategy, we propose to measure the degree of semantic informativity of deductions as a dynamic phenomenon, based on explicit definitions of insertions and deletions on databases (Section 2). As databases may not correspond to reality, this is the sense in which we do not assume truth as part of the \emph{definiens} of semantic information. More precisely, we do not use the concept of truth as correspondence in order to define semantic informativity (Cf. \cite{kirkham1992}). There is, even so, a sense in which truth is part of semantic information. Given the explicit definitions of structural operations, we will delimit a dynamic concept of informational complexity (Section 3), on the one hand, and a dynamic concept of informational relevancy (Section 4), on the other. Hence, semantic informativity will be conceived in terms of complexity and relevancy, which relies on the Tarskian conception of truth \cite{tarski1933}.

In other words, we sustain that truth is a constituent of semantic information, but it is a presupposition for informativity, not part of the definition of semantic information. In our approach false propositions can be semantically informative. We measure, however, the semantic informativity of propositions looking for what is sufficient to become them true. More accurately, we can say that semantic informativity is the product of two factors: how many structural operations we do for propositions to become true with respect to our database and how many of these propositions are relevant. As consequence of this approach, with respect to finite databases, we will exhibit a unified solution to Bar-Hillel-Carnap paradox and Hintikka's scandal of deduction (Section 5).

\section{Structural operations}

According to \cite{floridi2011}, semantic information is true well-formed data. We intend to analyze the semantic information of deductions. Since deductions are compounded of propositions, we can assume that data is expressed into propositions. Databases are, in turn, organized collections of data (Cf. \cite{kroenke2007}). Thus we can conceive that the semantic information of deductions is represented into databases. From a logical point of view, the concept of database could be understood by means of the concept of mathematical structure.

\bre

We write $A(\phi) = 1$ to indicate that proposition $\phi$ is true in structure $A$.

\ere

\bde\label{21}

Let $S$ be a signature. An \emph{$S$-database} is a pair $D=(A,T)$ compounded of a finite first-order structure $A$ over $S$ and a correct finite classical first-order theory $T$ over $S$ and about $A$. The theory $T$ is correct about $A$ in the sense that $A(\phi) = 1$ for all $\phi \in T$.

\ede

\bre

We use $\ulcorner\alpha\urcorner= \beta$ to express that the symbol $\beta$ is a formal representation of the expression $\alpha$. Besides, $X_{\beta}$ is the interpretation of $\beta$ in the structure $A$, where $X$ is a set over the domain of $A$. The domain of $A$ will be represented by $dom(A)$, and individuals of $dom(A)$ are represented by a bar above letters.

\ere

\bex\label{22}

Let $D_{0}=(A_{0},T)$ be a database with signature $S=(\{s,l,a\},\{C,E,H\})$, for $s=\ulcorner\mbox{S\~ao\ Paulo}\urcorner$, $l=\ulcorner\mbox{London}\urcorner$, $a=\ulcorner\mbox{Avenida\ Paulista}\urcorner$, $C=\ulcorner\mbox{City}\urcorner$,
$E=\ulcorner\mbox{Street}\urcorner$ and $H=\ulcorner\mbox{To\ have}\urcorner$, such that $A_{0}=(\{\bar{s},\bar{l},\bar{a}\},\bar{s}_{s},\bar{l}_{l},\bar{a}_{a},\{\bar{s},\bar{l}\}_{C},\{\bar{a}\}_{E},\{(\bar{s},\bar{a}),(\bar{l},\bar{a})\}_{H})$ and $T=\{\forall x(Cx\to\exists yHxy),\forall x(Cx\vee Ex),\neg El,Cs\}$.

\eex

The example \ref{22} indicates some features of the approach that we intend to develop here. The fact that $T$ is correct with
respect to $A$ does not exclude the non-correspondence of our database to reality. For instance, it is true in $A$ that $Hla\wedge Ea$, in words, it is true in $A$ that London has a city called Avenida Paulista, what is false until  the date of this paper. The theory $T$ represents the fundamental facts of our database that we consider true, that is to say, $T$ has our \emph{fundamental beliefs}. It is important to note that $T$ may not be complete about $A$, i.e., it is possible that not all true propositions about $A$ are in $T$; example \ref{22} shows this.

To measure semantic informativity of deductions, we define a dynamic of changes in databases through operations called \emph{structural operations}. The first one consists of introducing possibly new objects into the structure of the database and interpreting the possibly new symbols with these objects.

\bre

We presume that constants are the unique symbols of ariety $0$ and that relational and functional symbols have ariety greater than $0$.

\ere

\bde\label{23}

Let $D=(A,T)$ be an $S$-database. An \emph{insertion} of an $n$-ary symbol $\sigma$ into $D$ is a database $D'=(A,T)$ such that $A'$ is a structure over $S'=S\cup\{\sigma\}$ with the following properties: 

\begin{enumerate}

\item $A'(\tau)=A(\tau)$ for all $\tau\neq\sigma$ such that $\tau\in S$; 

\item If $n=0$ and $\sigma\in S$, then $dom(A') = dom(A)$ and $A'(\sigma) = A(\sigma)$, but if $n=0$ and $\sigma\notin S$, then $dom(A') = dom(A) \cup \{ a \}$ and $A'(\sigma)=a$, provided $A'(\phi)=1$ for all $\phi\in T$;

\item If $n>0$, then $dom(A') = dom(A) \cup \{a_{1},\ldots,a_{n}\}$ and $A'(\sigma)=A(\sigma)\cup\{(a_{1},\ldots,a_{n})\}$, provided $A'(\phi)=1$ for all $\phi\in T$. 

\end{enumerate}

\ede

\bex\label{24}

Let $D_{0}=(A,T)$ be the database in example \ref{22}. The database $D_{1}=(A_{1},T)$ with signature $S'=S\cup\{b\}$, where $b=\ulcorner\mbox{Shaftesbury\ Avenue}\urcorner$, and $A_{1}=(\{\bar{s},\bar{l},\bar{a}\},\bar{s}_{s},\bar{l}_{l},\bar{a}_{a},\bar{a}_{b},\{\bar{s},\bar{l}\}_{C},\{\bar{a}\}_{E},\{(\bar{s},\bar{a}),(\bar{l},\bar{a})\}_{H})$
is an insertion of $b$ into $D$. On the other hand, given the $S'$-structure $A_{2}=(\{\bar{s},\bar{l},\bar{a},\bar{b}\},\bar{s}_{s},\bar{l}_{l},\bar{a}_{a},\bar{a}_{b},\{\bar{s},\bar{l}\}_{C},\{\bar{a},\bar{b}\}_{E},\{(\bar{s},\bar{a}),(\bar{l},\bar{a})\}_{H})$, $D_{2}=(A_{2},T)$ is an insertion of $E$ into $D_{1}$. Nonetheless, considering the $S'$-structure $A^{*}=(\{\bar{s},\bar{l},\bar{a},\bar{b}\},\bar{s}_{s},\bar{l}_{l},\bar{a}_{a},\bar{b}_{b},\{\bar{s},\bar{l}\}_{C},\{\bar{a}\}_{E},\{(\bar{s},\bar{a}),(\bar{l},\bar{a})\}_{H})$, we have that $D^{*}=(A^{*},T)$ is not an insertion of $b$ into $D_{1}$ because $A^{*}(\forall x(Cx\vee Ex)) = 0$.

\eex

The example \ref{24} shows that it is not necessary the introduction of a new object into the structure of the database in order to define an insertion (Cf. database $D_{1}$). It is sufficient to add an element in interpretation of some symbol. On the other hand, this example also exhibits that the introduction of a new object into the structure of the database does not necessarily produce an insertion (Cf. database $D^{*}$). It is necessary to guarantee that the fundamental beliefs of the database are true in the new structure.

We turn now to the second structural operation. This operation consists of removing possibly old objects from the structure of the database and interpreting the possibly new symbols with the remain objects.

\bde

Let $A$ be an $S$-structure. An element $a \in dom(A)$ is called \emph{free} for the symbol $\sigma\in S$ if there is no
constant $\tau\in S$ such that $\tau\neq\sigma$ and $A(\tau)=a$ or symbol $\delta$ with arity $n>0$ such that $(b_{1},\ldots,b_{n})\in A(\delta)$ and $a=b_{i}$ for some $1 \leq i \leq n$. If $\sigma$ is an $n$-ary symbol, we write $A(\sigma)_i$ to denote element $a_i$ of some sequence $(a_1, \ldots, a_n) \in A(\sigma)$.

\ede

\bde\label{25}

Let $D=(A,T)$ be an $S$-database. A \emph{deletion} of an $n$-ary symbol $\sigma$
from $D$ is a database $D'=(A',T)$ such that $A'$ is a structure over $S'$, where either $S' = S$ or $S' = S - \{\sigma\}$, with the following properties:

\begin{enumerate}

\item $A'(\tau) = A(\tau)$ for all $\tau\neq\sigma$ such that $\tau\in S-\{\sigma\}$; 

\item If $n=0$ and $S'=S$, then $dom(A') = dom(A)$ and $A'(\sigma)=a$ for some $a \in dom(A)$ such that $a\neq A(\sigma)$, but if $n=0$ and $S' = S - \{\sigma\}$, then $dom(A') = dom(A) -\{A(\sigma)\}$ if $A(\sigma)$ is free for $\sigma$, otherwise, $dom(A') = dom(A)$, in any case provided $A'(\phi)=1$ for all $\phi\in T$; 

\item If $n>0$ and $S'=S$, then $dom(A') = dom(A)$ and $A'(\sigma) = A(\sigma) - \{(a_{1}, \ldots,a_{n})\}$, but if $n>0$ and $S' = S - \{\sigma\}$, then $dom(A') = dom(A) - \{A(\sigma)_{1}, \ldots, A(\sigma)_{k}\}$, where each $A(\sigma)_{i}$ is free for $\sigma$ and , in any case provided $A'(\phi) = 1$ for all $\phi \in T$.

\end{enumerate}

\ede

\bex\label{26}

Let $D_{0}=(A,T)$ be the database in example \ref{22}. The database $D'_{1}=(A'_{1},T)$ with signature $S$ and $A'_{1}=(\{\bar{s},\bar{l},\bar{a}\},\bar{a}_{s},\bar{l}_{l},\bar{a}_{a},\{\bar{s},\bar{l}\}_{C},\{\bar{a}\}_{E},\{(\bar{s},\bar{a}),(\bar{l},\bar{a})\}_{H})$ is a deletion of $s$ from $D$. On the other hand, $D'_{2}=(A'_{2},T)$ is a deletion of $H$ from $D'_{1}$ where $A'_{2}=(\{\bar{s},\bar{l},\bar{a}\},\bar{a}_{s},\bar{l}_{l},\bar{a}_{a},\{\bar{s},\bar{l}\}_{C},\{\bar{a}\}_{E},\{(\bar{l},\bar{a})\}_{H})$ is an $S$-structure. Nonetheless, for $A'_{3}=(\{\bar{l},\bar{a}\},\bar{a}_{s},\bar{l}_{l},\bar{a}_{a},\{\bar{l}\}_{C},\{\bar{a}\}_{E},\{(\bar{l},\bar{a})\}_{H})$, a structure over the signature $S$, we have that $D'_{3}=(A_{3},T)$ is not a deletion of $C$ from $D'_{1}$ because in this case, in despite of $A'_{3}(\phi) = 1$ for $\phi\in T$, $D'_{3}(H) \neq D'_{1}(H)$. Note, however, that $D'_{3}$ is a deletion of $C$ from $D'_{2}$.

\eex

The example \ref{26} illustrates other important aspects of databases. First, each structural operation can change the interpretation of one, and only one, symbol. For this reason $D'_{3}$ is not a deletion of $C$ from $D'_{1}$, but $D'_{3}$ is a deletion of $C$ from $D'_{2}$. Second, we can only delete the elements of the symbol's interpretation to be deleted, without interfering into the interpretations of the others. The restriction to free objects ensures that.

Insertions and deletions on databases are well-known primitive operations in the study of databases and dynamic computational complexity (Cf. \cite{kroenke2007}, \cite{weber2007}). Nevertheless, to the best of our knowledge, they have being conceived as undefined notions. In contrast, we presented a logical conception about databases with explicit definitions of insertion and deletion. These definitions will permit us to develop a quantitative approach to semantic informativity in the next sections.

\section{Informational complexity}

For Popper, ``the amount of empirical information conveyed by a theory, or its \emph{empirical content}, increases with its degree of falsifiability'' \cite{popper1959}[pg.96]. Thus, by duality we can infer that the empirical information of a theory decreases with its degree of verifiability. It is well-known that Popper does not assign meaning to the notion of \emph{degree of verifiability}, but in the context of our approach we can do that. Indeed, the number of structural operations performed on our database $D = (A,T)$ to produce a database $D=(B,T)$ in which $\phi$ is true can be viewed as the degree of verifiability of $\phi$. In this case, the more structural operations $\phi$ requires to become true the less verifiable $\phi$ is. In what follows we provide a precise meaning to that perspective.

\bde\label{31}

An \emph{update} of an $S$-database $D$ is an at most enumerable sequence $\mathcal{D}=(D_{i}:i<\alpha)$ such that $D_{0}=D$ and each $D_{i+1}$ is an insertion or deletion in $D_{i}$. An update $\mathcal{D}=(D_{i}:i<\alpha)$ of $D$ is \emph{satisfactory} for a proposition $\phi$ if $\alpha=n+1$, for some $n \in \mathbb{N}$, $D_{n+1}=(A_{n},T)$ and $A_{n}(\phi)=1$; otherwise, $\mathcal{D}$ is said to be \emph{unsatisfactory} for $\phi$. We write $\mathcal{D}(\phi)=1$ or $\mathcal{D}(\phi)=0$ according to $\mathcal{D}$ is satisfactory for $\phi$ or not. A proposition $\phi$ is \emph{$D$-acceptable} if there is an update of $D$ satisfactory for $\phi$, otherwise, $\phi$ is called \emph{$D$-unacceptable}. Besides that, if $\mathcal{D} = ( D_i : i < n )$, with $D_{i}=(A_{i},T)$, and $\mathcal{D}(\phi) = 1$, we define $\| \mathcal{D} \| = \min \{ i < n : A_i(\phi) = 1 \}$.

\ede

\bex\label{32}

Let $D_{0}$ be the database in example \ref{22} and $D_{1}$ in example \ref{24}. The sequence $\mathcal{D}=(D_{0},D{}_{1})$ is an update of $D$ satisfactory for $Eb$ and $Hlb$. Consider the databases $D'_{1}$, $D'_{2}$ and $D'_{3}$ in example \ref{26}. The sequence $\mathcal{D}'=(D_{0},D'_{1},D'_{2},D'_{3})$ is an update of $D$ satisfactory for $Es\wedge\neg Hsa$ but not for $s=a$ because the last proposition is false in $A'_{3}=(\{\bar{l},\bar{a}\},\bar{a}_{s},\bar{l}_{l},\bar{a}_{a},\{\bar{l}\}_{C},\{\bar{a}\}_{E},\{(\bar{l},\bar{a})\}_{H})$.

\eex

In other words, an update satisfactory for a proposition $\phi$ is a sequence of changes made on our database $D$ to produce a structure $B$ in which $\phi$ is true. This number of chances measures a certain complexity of $\phi$ and its degree of verifiability is inverse to that complexity. For the sake of generality, in the definition below we formalize this notion of complexity with respect to sets of updates (the case of one update is so a particular sub-case).

\bde\label{33}

Let $D$ be a database and let $\bar{\mathcal{D}}$ be a collection of updates over $D$. The \emph{informational complexity} of $\phi$ in $\bar{\mathcal{D}}$ is defined by 

\[
C_{\bar{\mathcal{D}}}(\phi)=\min_{\mathcal{D}\in\bar{\mathcal{D}}}\{ \| \mathcal{D} \|:\mathcal{D}(\phi) = 1\},
\]

if there is $\mathcal{\mathcal{D}}\in\bar{\mathcal{D}}$ satisfactory for $\phi$, otherwise, we define that

\[
C_{\bar{\mathcal{D}}}(\phi)=0.
\]

\ede

\bex\label{34}

Let $\bar{\mathcal{D}}=\{\mathcal{D}\}$, where $\mathcal{D}$ is the update in example \ref{32}. Since $\mathcal{D}$ is satisfactory for $Eb$ and $Hlb$, we have $C_{\bar{\mathcal{D}}}(Eb)=C_{\bar{\mathcal{D}}}(Hlb)=1$ and so $C_{\bar{\mathcal{D}}}(Eb\wedge Hlb)=C_{\bar{\mathcal{D}}}(Eb\vee Hlb)=1$. Now, let $\bar{\mathcal{D}}'=\{\mathcal{D}'\}$, where $\mathcal{D}'$
is the update in example \ref{32}. Then, we have $C_{\bar{\mathcal{D}}'}(\neg Cs) = 1$, $C_{\bar{\mathcal{D}}'}(\neg Hsa) = 3$, $C_{\bar{\mathcal{D}}'}(s=a)=0$ and so $C_{\bar{\mathcal{D}}'}(\neg Cs\wedge\neg Hsa)=3$ but $C_{\bar{\mathcal{D}}'}(\neg Cs\wedge s=a)=0$.

\eex

The example \ref{34} exhibits that, given an update, different propositions can have different complexities, but different propositions can have same complexity too. In special, an intriguing point deserves attention. It seems natural to think that $Eb\wedge Hlb$ is in some sense more complex than $Eb$ and $Hlb$. Here we do not have this phenomenon. Since $C_{\bar{\mathcal{D}}}(Eb)=C_{\bar{\mathcal{D}}}(Hlb) = C_{\bar{\mathcal{D}}}(Eb\wedge Hlb)=C_{\bar{\mathcal{D}}}(Eb\vee Hlb)=1$, informational complexity \emph{is not} an additive measure of complexity of propositions. That is to say, informational complexity of conjunctions or disjunctions is not the sum of their components. Another interesting point is that $C_{\bar{\mathcal{D}}'}(\neg Hsa)>C_{\bar{\mathcal{D}}'}(\neg Cs)$ but $C_{\bar{\mathcal{D}}'}(\neg Hsa)=C_{\bar{\mathcal{D}}}(\neg Cs\wedge Hlb)=3$. This reflects the fact that updates are sequences. First, we took $\neg Cs$ and made it $D$-acceptable, later we made $\neg Hsa$ an $D$-acceptable proposition. When $\neg Hsa$ is $D$-satisfactory there is nothing more to be done, as far as the conjunction $\neg Cs\wedge\neg Hsa$ is concerned\footnote{In a future work, we will analyze these characteristics of the notion of complexity in details.}.

We can now define the informational complexity of sets of propositions with respect to a given set of updates over a database.

\bde

Let $D$ be a database and let $\bar{\mathcal{D}}$ be a collection of updates over $D$. The \emph{informational complexity} of a set of formulas $\{\phi_{0},\ldots,\phi_{n}\}$ in $\bar{\mathcal{D}}$ is defined by

\[
C_{\bar{\mathcal{D}}}(\{\phi_{0},\ldots,\phi_{n}\})=\sum_{i=0}^{n}C_{\bar{\mathcal{D}}}(\phi_{i}).
\]

We define the informational complexity of a deduction $\phi_{0},\ldots,\phi_{n}$ as the informational complexity of the set $\{\phi_{0},\ldots,\phi_{n}\}$.

\ede

From these definitions we get an important result for our analysis of Bar-Hillel-Carnap paradox and Hintikka's scandal of deduction.

\bla\label{35}

For every collection of updates $\bar{\mathcal{D}}$ over an $S$-database $D$, the following statements are true:

\begin{enumerate}

	\item If $\phi$ is a tautology over $S$ and $\bar{\mathcal{D}}$ is not empty, then $C_{\bar{\mathcal{D}}}(\phi) = 0$.

	\item If $\phi$ is a tautology with symbols not in $S$ and $\bar{\mathcal{D}}$ has some satisfactory update for $\phi$, then $C_{\bar{\mathcal{D}}}(\phi) > 0$.

	\item if $\phi$ is $D$-unacceptable, then $C_{\bar{\mathcal{D}}}(\phi) = 0$.

\end{enumerate}

\ela

\bpr

Let $D=(A,T)$ be an $S$-database.  If $\bar{\mathcal{D}}$ is not empty, it is immediate from the definition of $C_{\bar{\mathcal{D}}}(\phi)$ that $C_{\bar{\mathcal{D}}}(\phi) = 0$ for $\phi$ a tautology over $S$, because $A(\phi)=1$ and every update $\mathcal{D}=(D_{i}:i < \alpha)$ over $D$ is such that $D_{0} = (A,T)$.

Suppose that $\phi$ is a tautology not over $S$. Let $\bar{\mathcal{D}}$ be a collection of updates over $D$ with one satisfactory update for $\phi$. Then, there are symbols in $\phi$ that are not in $S$, say that these symbols are $\sigma_{0},\ldots,\sigma_{k}$. In this case, every structure $B$ such that $B\vDash\psi$ should be a structure over a signature $S'\supseteq S\cup\{\sigma_{0},\ldots,\sigma_{k}\}$. This means that $B$ has interpretations that $A$ does not have, namely $B(\sigma_{0}),\ldots,B(\sigma_{k})$. Therefore, every update $\mathcal{D}=(D_{0},\ldots,D_{n})$ in $\bar{\mathcal{D}}$ satisfactory for $\phi$ is such that $D_{0}=D$ and $D_{n}=(B,T)$ should have $n>k$, for it is necessary at least one structural operation to define each $B(\sigma_{i})$. This shows that $0<C_{\bar{\mathcal{D}}}(\phi)$.

If $\phi$ is $D$-unacceptable, there is no update $\mathcal{D}=(D_{0},\ldots,D_{n})$ of $D$ satisfactory for $\phi$. This means that $C_{\bar{\mathcal{D}}}(\phi)=0$ for every collection of updates $\bar{\mathcal{D}}$ over $D$.

\epr

We have defined $C_{\bar{\mathcal{D}}}(\phi)$ in such a way that, in principle, if $\phi$ was a tautology or contradiction, then $C_{\bar{\mathcal{D}}}(\phi)=0$. Nonetheless, lemma \ref{35} shows that, from an informational point of view, there is a natural asymmetry between tautologies and contradictions. Tautologies could be complex in some situations, notably, when they have new symbols (relatively to a given database). In contrast, contradictions always have null informational complexity.

\section{Informational relevancy}

The fundamental distinction between a set of propositions and a deduction is that in a deduction there is a necessary connection among its propositions and one special proposition called \emph{conclusion} (Cf. \cite{shapiro2005}). In the process of obtaining an update satisfactory for the conclusion of a deduction, some propositions may be relevant with respect to our beliefs, but others not. The informativity of a deduction depends on the relevance we associated to its propositions. For this reason, in what follows we will define a logical notion of relevance.

\bde\label{41}

Let $\mathcal{D}$ be an update of a database $D=(A,T)$. The \emph{informational relevant propositions} of a set of formulas $\{\phi_{0},\ldots,\phi_{n}\}$ are those for which $\mathcal{D}$ is satisfactory but that are not logical consequences of $T$, i.e., the propositions in the set

\[
\mathcal{D}(\{\phi_{0},\ldots,\phi_{n}\})=\{\phi\in\{\phi_{0},\ldots,\phi_{n}\}:\mathcal{D}(\phi)=1\mbox{ and }T\nvDash\phi\}.
\]

\ede

From now on we will omit the parentheses in $\mathcal{D}(\{\phi_{0},\ldots,\phi_{n}\})$.

\bex\label{42}

For $\mathcal{D}=(D_{0})$, where $D_{0}$ is the database in example \ref{32}, we have $\mathcal{D}(Ea,\exists xEx)=\{Ea\}$. Let $\mathcal{D}'=(D_{0},D_{1})$, where $D_{0}$ and $D_{1}$ are the databases in example \ref{32} too. In this case, $\mathcal{D}'(\forall x(Cx\to\neg Ex),Cb,Cb\to\neg Eb,\neg Eb)=\{\forall x(Cx\to\neg Ex),Cb\to\neg Eb\}$.

\eex

In definition \ref{41} we have opted by a strong requirement about relevancy of propositions, namely, only the non-logical consequences of our fundamental believes are relevant. For practical purposes, this requirement could be relaxed. Here we are interested in foundational problems, Bar-Hillel-Carnap paradox and Hintikka's scandal of deduction. In this case definition \ref{41} seems to be a good beginning, but there is a sense in which it is too much weak for analyzing information of deductions.

\bpo\label{43}

For every update $\mathcal{D}$ of $D$ and valid deduction $\phi_{0},\ldots,\phi_{n}$, if $\mathcal{D}$ is unsatisfactory for $\phi$, then some premisses of $\phi_{0},\ldots,\phi_{n}$ are not informational relevant.

\epo

\bpr

Suppose that $\Gamma=\{\phi_{0},\ldots,\phi_{k}\}$ are the premisses of $\phi_{0},\ldots,\phi_{n}$, and consider an update $\mathcal{D}$ of $D$ unsatisfactory for its conclusion $\phi_{n}$. If $\phi_{0},\ldots,\phi_{n}$ is a valid deduction, then, for every $A$, if $A(\phi)=1$ for every $\phi\in\Gamma$, then $A(\phi_{n})=1$. As $\mathcal{D}(\phi_{n})=0$, it follows that $\mathcal{D}(\phi)=0$ for some $\phi\in\Gamma$. Hence, $\mathcal{D}(\Gamma)<|\Gamma|$, where $|\Gamma|$ indicates the cardinality of $\Gamma$.

\epr

Hence, if we analyze which propositions of a valid deduction are relevant in an update unsatisfactory for its conclusion, we will establish \emph{a priori} that some premiss is irrelevant for the deduction. In other words, we loose the possible information associated to the connection between the conclusion and some premiss of the deduction under analysis. This is untenable. We should, then, analyze the relevancy of deductions with respect to updates satisfactory for conclusions. Moreover, if $\phi_{0},\ldots,\phi_{n}$ is an invalid deduction, then there is always an update $\mathcal{D}$ unsatisfactory for its conclusion $\phi_{n}$, and so invalid deductions are not informative, for they do not guarantee any intrinsic relationship among conclusion and premisses.

Let us refine a little more the remark above. Given a deduction $\phi_{0},\phi_{1},\ldots,\phi_{n}$ compounded of propositions over a signature $S$, we assume that its premisses are the first propositions $\phi_{0},\phi_{1},\ldots,\phi_{k}$ ($k\leq n$) and the conclusion is the last one $\phi_{n}$. We call this set the \emph{support} of $\phi_{0},\phi_{1},\ldots,\phi_{n}$, that is to say, the support of a deduction is its set of premisses plus its conclusion, in symbols $\{\phi_{0},\ldots,\phi_{k}\}\{\phi_{n}\}$. Now suppose that we have a valid deduction $\phi_{0},\phi_{1},\ldots,\phi_{n}$ with support $\{\phi_{0},\ldots,\phi_{k}\}\{\phi_{n}\}$. Then, as $\phi_{0},\phi_{1},\ldots,\phi_{n}$ is valid, $\mathcal{D}(\phi_{i})=1$, $0\leq i\leq k$, implies that $\mathcal{D}(\phi_{j})=1$, $k\leq j\leq n$, for every update $\mathcal{D}$. Hence, although the steps from the premisses to the conclusion of $\phi_{0},\phi_{1},\ldots,\phi_{n}$ may be relevant from a syntactical standpoint, from a semantic perspective we can only pay attention to the relationship among premisses of $\phi_{0},\phi_{1},\ldots,\phi_{n}$ and its conclusion, that is to say, we can just look at the support $\{\phi_{0},\ldots,\phi_{k}\}\{\phi_{n}\}$. 

\bde\label{44}

Let $D$ be a database, $\bar{\mathcal{D}}$ a collection of updates over $\mathcal{D}$ and $\phi_{0},\ldots,\phi_{n}$ a deduction with support $\{\phi_{0},\ldots,\phi_{k}\}\{\phi_{n}\}$. If $\phi_{0},\ldots,\phi_{n}$ is valid, then the \emph{informational relevancy} of $\phi_{0},\ldots,\phi_{n}$ in $\bar{\mathcal{D}}$ is defined by 

\[
R_{\bar{\mathcal{D}}}(\phi_{0},\ldots,\phi_{n})=\frac{|\mathcal{D}(\{\phi_{0},\ldots,\phi_{k}\}\{\phi_{n}\})|}{|\{\phi_{0},\ldots,\phi_{k}\}\{\phi_{n}\}|},
\]

where $\mathcal{D}$ is the smallest update in $\bar{\mathcal{D}}$ satisfactory for $\phi_{n}$, but if there is no such a $\mathcal{D}$ or $\phi_{0},\ldots,\phi_{n}$ is invalid, then we define that 

\[
R_{\bar{\mathcal{D}}}(\phi_{0},\ldots,\phi_{n})=0.
\]

\ede

\bex\label{45}

Let $\bar{\mathcal{D}}=\{\mathcal{D},\mathcal{D}'\}$, where $\mathcal{D}$ and $\mathcal{D}'$ are the updates in example \ref{42}. Hence, $R_{\bar{\mathcal{D}}}(Ea,\exists xEx)=1$, but $R_{\bar{\mathcal{D}}}(\forall x(Cx\to\neg Ex),Cb,Cb\to\neg Eb,\neg Eb)=0$ because there is no update in $\bar{\mathcal{D}}$ satisfactory for $\neg Eb$. Nonetheless, consider the new update $\mathcal{D}''=(D_{0},D_{1},D_{2},D_{3},D_{4})$ such that $D_{0}$, $D_{1}$ and $D_{2}$ are the updates in example
\ref{24}, $A_{3}=(\{\bar{s},\bar{l},\bar{a},\bar{b}\},\bar{s}_{s},\bar{l}_{l},\bar{a}_{a},\bar{b}_{b},\{\bar{s},\bar{l}\}_{C},\{\bar{a},\bar{b}\}_{E},\{(\bar{s},\bar{a}),(\bar{l},\bar{a})\}_{H})$ and $A_{4}=(\{\bar{s},\bar{l},\bar{a},\bar{b}\},\bar{s}_{s},\bar{l}_{l},\bar{a}_{a},\bar{b}_{b},\{\bar{s},\bar{l}\}_{C},\{\bar{a}\}_{E},\{(\bar{s},\bar{a}),(\bar{l},\bar{a})\}_{H})$. Now, put $\bar{\mathcal{D}}'=\{\mathcal{D},\mathcal{D}',\mathcal{D}''\}$. In this case, we have that $R_{\bar{\mathcal{D}}'}(\forall x(Cx\to\neg Ex),Cb,Cb\to\neg Eb,\neg Eb)=2/3$ because $\mathcal{D}''(\{\forall x(Cx\to\neg Ex),Cb\}\{\neg Eb\})=\{\forall x(Cx\to\neg Ex),\neg Eb\}$.

\eex

In example \ref{45} the fact that $R_{\bar{\mathcal{D}}}(\forall x(Cx\to\neg Ex),Cb,Cb\to\neg Eb,\neg Eb) = 0$ show us that we can have valid deductions with new symbols but irrelevant. On the other hand, $R_{\bar{\mathcal{D}}}(Ea,\exists xEx) = 1$ shows that it is not necessary to consider new symbols to find deductions with non-null relevancy. This means that informational relevancy has a criteria of nullity different from informational complexity.

\bla\label{46}

For every set of updates $\bar{\mathcal{D}}$ over an $S$-database $D=(A,T)$, if $\phi_{0},\ldots,\phi_{n}$ is a deduction with support $\{\phi_{0},\ldots,\phi_{k}\}\{\phi_{n}\}$, then $R_{\bar{\mathcal{D}}}(\phi_{0},\ldots,\phi_{n})=0$ in the following cases:

\begin{enumerate}

\item $T$ is a complete theory of $A$ and $\phi_{0},\ldots,\phi_{n}$ is an $S$-deduction;

\item $\phi_{n}$ is a tautology and $\{\phi_{0},\ldots,\phi_{k}\}=\oslash$;

\item $\phi_{n}$ is $D$-unacceptable.

\end{enumerate}

\ela

\bpr

Let $\phi_{0},\ldots,\phi_{n}$ be a deduction with support $\{\phi_{0},\ldots,\phi_{k}\}\{\phi_{n}\}$. If $\phi_{0},\ldots,\phi_{n}$ is invalid, then, by definition, $R_{\bar{\mathcal{D}}}(\phi_{0},\ldots,\phi_{n})=0$ for every $\bar{\mathcal{D}}$. Besides that, if $\bar{\mathcal{D}}$ is empty, by definition $R_{\bar{\mathcal{D}}}(\phi_{0},\ldots,\phi_{n})=0$. Then, let us consider that $\phi_{0},\ldots,\phi_{n}$ is valid and $\bar{\mathcal{D}}$ is not empty.

Suppose that $T$ is a complete theory of $A$ and $\phi_{0},\ldots,\phi_{n}$ is an $S$-deduction. Then, for every $\phi\in\{\phi_{0},\ldots,\phi_{k}\}\{\phi_{n}\}$, if $B$ is an $S$-structure that is a model of $T$ and $B(\phi)=1$,
then $T\vDash\phi$. By definition, every update of $D$ is compounded of models of $T$. Thus, if $\mathcal{D}=(D_{0},\ldots,D_{m})$ is an update of $D$ and $A_{m}(\phi)=1$, then $T\vDash\phi$ and so $\phi$ is not relevant, but if $A_{m}(\phi)=0$, then $\phi$ is not relevant too. Therefore, if $T$ is a complete theory of $A$ and $\phi_{0},\ldots,\phi_{n}$ is an $S$-deduction, then $R_{\bar{\mathcal{D}}}(\phi_{0},\ldots,\phi_{n})=0$ for every set of updates $\bar{\mathcal{D}}$ over the $S$-database
$D=(A,T)$.

Suppose $\phi_{n}$ is a tautology and $\{\phi_{0},\ldots,\phi_{k}\}=\oslash$. If $\phi_{n}$ is a tautology over the signature $S$, then $T\vDash\phi_{n}$ and, in virtue of $\{\phi_{0},\ldots,\phi_{k}\}=\oslash$, $R_{\bar{\mathcal{D}}}(\phi_{0},\ldots,\phi_{n})=0$ for every set of updates $\bar{\mathcal{D}}$ over the $S$-database $D=(A,T)$. If $\phi_{n}$ is a tautology with at least one symbol not in $S$, say that $S'$ is the signature of $\phi_{n}$, then we cannot affirm that $T\vDash\phi_{n}$ because $\phi_{n}$ is not in the language of $T$. Nevertheless, it is impossible to assert $T\nvDash\phi_{n}$.
Indeed, suppose that $T\nvDash\phi_{n}$. Hence, there is an $S\cup S'$-structure $B$ such that $B$ is model of $T$ and it is not a model of $\phi_{n}$, but this an absurd: $\phi_{n}$ is a tautology over $S\cup S'$ and for that reason $B(\phi_{n})=1$. Therefore, if $\bar{\mathcal{D}}$ is a set of updates over the $S$-database $D=(A,T)$ with some $\mathcal{D}$
satisfactory for $\phi_{n}$, then $R_{\bar{\mathcal{D}}}(\phi_{0},\ldots,\phi_{n})=0$.

If is $D$-unacceptable, then there is not update $\mathcal{D}$ of $D$ such that $\mathcal{D}(\phi_{n})=1$ and, by definition, $R_{\bar{\mathcal{D}}}(\phi_{0},\ldots,\phi_{n})=0$ for every set of updates $\bar{\mathcal{D}}$ over the $S$-database
$D=(A,T)$.

\epr

This result shows that deductions can be relevant only when we do not have a complete theory of the structure of the database. Moreover, \emph{qua} deductions, isolated logical facts (tautologies and contradictions) have no relevance. Hence, informational relevancy is an inferential notion. This conclusion will be crucial in our solution to Bar-Hillel-Carnap paradox and Hintikka's scandal of deduction.

\section{Semantic informativity}

Having at hand dynamic concepts of complexity and relevance, it seems reasonable to state two statements. On the one hand, the more complex the propositions in support of a deduction are the more informative the deduction is. On the other hand, the more relevant the support of a deduction is the more information it provides. Our definition of semantic informativity of deductions relies on that idea.

\bde\label{51}

The \emph{semantic informativity} $I_{\bar{\mathcal{D}}}(\phi_{0},\ldots,\phi_{n})$ of a deduction $\phi_{0},\ldots,\phi_{n}$ in the set of updates $\bar{\mathcal{D}}$ over the database $D$ is defined by 

\[
I_{\bar{\mathcal{D}}}(\phi_{0},\ldots,\phi_{n})=C_{\bar{\mathcal{D}}}(\phi_{0},\ldots,\phi_{n})R_{\bar{\mathcal{D}}}(\phi_{0},\ldots,\phi_{n}).
\]

\ede

\bex\label{52}

Let $\bar{\mathcal{D}}=\{\mathcal{D},\mathcal{D}',\mathcal{D}''\}$, where $\mathcal{D}$, $\mathcal{D}'$ and $\mathcal{D}''$ are the updates in example \ref{45}. Then, $I_{\bar{\mathcal{D}}}(Ea,\exists xEx)=1\cdot1=1$ and $I_{\bar{\mathcal{D}}'}(\forall x(Cx\to\neg Ex),Cb,Cb\to\neg Eb,\neg Eb)=(2/3)\cdot4=8/3$.

\eex

In other words, the semantic informativity of a deduction $\phi_{0},\ldots,\phi_{n}$ is directly proportional to the complexity and relevance of its propositions. This standpoint permit us to conceive the semantic informativity of a proposition as a special case of the semantic informativity of deductions. What is more, we can prove that the semantic informativity of a proposition is its informational complexity.

\bte\label{53}

For every collection of updates $\bar{\mathcal{D}}$ over an $S$-database $D$, $I_{\bar{\mathcal{D}}}(\phi)=C_{\bar{\mathcal{D}}}(\phi)$.

\ete

\bpr

Let $\bar{\mathcal{D}}$ be a collection of updates over an $S$-database $D=(A,T)$. By definition, $C_{\bar{\mathcal{D}}}(\{\phi\})=C_{\bar{\mathcal{D}}}(\phi)$. If $C_{\bar{\mathcal{D}}}(\phi)=0$, then $I_{\bar{\mathcal{D}}}(\phi)=R_{\bar{\mathcal{D}}}(\phi)C_{\bar{\mathcal{D}}}(\phi)=0$. If $C_{\bar{\mathcal{D}}}(\phi)>0$, then there is an update $\mathcal{D}\in\bar{\mathcal{D}}$ such that $\mathcal{D}(\phi)=1$ but $A(\phi)=0$. Let us suppose that this $\mathcal{D}$ is the least one for which $\mathcal{D}(\phi)=1$. In this case, as $A$ is a model of $T$, it should be that $T\nvDash\phi$.
Therefore, $R_{\bar{\mathcal{D}}}(\{\phi\})=1$, which means that $I_{\bar{\mathcal{D}}}(\phi)=C_{\bar{\mathcal{D}}}(\phi)$.

\epr

\bex\label{54}

Let $\bar{\mathcal{D}}=\{\mathcal{D},\mathcal{D}',\mathcal{D}''\}$, where $\mathcal{D}$, $\mathcal{D}'$ and $\mathcal{D}''$ are the updates in example \ref{45}. We have that $I_{\bar{\mathcal{D}}}(Ea)=I_{\bar{\mathcal{D}}}(\exists xEx)=I_{\bar{\mathcal{D}}}(Ea\to\exists xEx)=0$, despite the fact that $I_{\bar{\mathcal{D}}}(Ea,\exists xEx)=1$.
Similarly, we also have that $I_{\bar{\mathcal{D}}}(\forall x(Cx\to\neg Ex))=0$, $I_{\bar{\mathcal{D}}}(Cb)=0$, $I_{\bar{\mathcal{D}}}(\neg Eb)=4$, $I_{\bar{\mathcal{D}}}((\forall x(Cx\to\neg Ex)\wedge Cb\wedge Cb\to\neg Eb)\to\neg Eb)=0$, although $I_{\bar{\mathcal{D}}}(\forall x(Cx\to\neg Ex),Cb,Cb\to\neg Eb,\neg Eb)=8/3$.

\eex

Theorem \ref{53} shows that semantic informativity of propositions is a measure of semantic information in Floridi's sense \cite{floridi2004}. It measures how many structural operations we do in order to obtain the semantic information of a proposition. For this reason, in example \ref{54} we have that $I_{\bar{\mathcal{D}}}(Cb)=0$: false well-formed data is not semantically informative; it should be true in our database. In contrast to Floridi, in our approach not all true well-defined data is informative. In example \ref{54}, $I_{\bar{\mathcal{D}}}(\forall x(Cx\to\neg Ex))=0$, but $\forall x(Cx\to\neg Ex)$ is true in the original update $D$. Furthermore, the realistic or anti-realistic nature of semantic information is out of question, for the propositions in our database may not be true in world. Differently from Floridi's theory, our theory of semantic informativity is strictly logical.

Since the notion of semantic informativity that we have proposed is logical, there is one more important distinction with respect to other recent theories of semantic information - not only with respect to Floridi's theory but also to \cite{sequoiah2009,jago2009,dalfonso2011}. Given a deduction $\phi_{0},\ldots,\phi_{n}$ and a set of updates $\bar{\mathcal{D}}$ over $D$, if we have $C_{\bar{\mathcal{D}}}(\phi_{0},\ldots,\phi_{n})=0$ or $R_{\bar{\mathcal{D}}}(\phi_{0},\ldots,\phi_{n})=0$, then the semantic informativity of $\phi_{0},\ldots,\phi_{n}$ is zero, it does not matter how $\phi_{0},\ldots,\phi_{n}$ is. Now, if $C_{\bar{\mathcal{D}}}(\phi_{n})=0$, then, by definition, $R_{\bar{\mathcal{D}}}(\phi_{0},\ldots,\phi_{n})=0$. In other words, if our database is $D$ and we intend to evaluate $I_{\bar{\mathcal{D}}}(\phi_{0},\ldots,\phi_{n})$ in some set of updates $\bar{\mathcal{D}}$ over $D$, we should look for a $\mathcal{D}$ in $\bar{\mathcal{D}}$ satisfactory for $\phi_{n}$, i.e., a $\mathcal{D}$ for which $C_{\mathcal{D}}(\phi_{n})>0$. Therefore, due to our notion of informational relevancy, the analysis of semantic informativity is also oriented to the conclusion of deductions. Theorem \ref{53} is a consequence of that inferential orientation.

In \cite{floridi2008}, Floridi already analyzed the relation between epistemic relevancy and semantic information, but he did not conceive relevancy a component of semantic informativity, as we have done here (Cf. \cite{allo2014}). In \cite{mares2009}, it was proposed to study informativity in relevant logics, but our concept of relevancy is not connected with non-classical logics. It is actually very far from the one used in relevant logics, for us a proposition can be relevant despite it is not used to obtain a conclusion - what is quite opposite to relevant logics. Our theory is a hegelian conception about semantic information. It analyzes semantic informativity in semantic-inferential terms, like Brandom \cite{brandom1989} did with respect to linguistic meaning. Whereas Floridi's conception is kantian, in the sense that it analyzes the relationship between propositions and the world in order to understand the conditions for semantic information.

Given that changes in the strong conception about semantic information, we can, first, solve Bar-Hillel-Carnap paradox.

\bte\label{55}

For some collections of updates $\bar{\mathcal{D}}$ over an $S$-database $D$, there are tautologies $\phi$ such that $I_{\bar{\mathcal{D}}}(\phi) = 0$ but there are tautologies $\psi$ such that  $I_{\bar{\mathcal{D}}}(\psi) > 0$. In contrast, for every collection of updates $\bar{\mathcal{D}}$ over an $S$-database $D$ and every contradiction $\phi$, $I_{\bar{\mathcal{D}}}(\phi) = 0$.

\ete
\bpr
By theorem \ref{53}, $I_{\bar{\mathcal{D}}}(\phi) = C_{\bar{\mathcal{D}}}(\phi)$ for every collection of updates $\bar{\mathcal{D}}$ over an $S$-database $D$. Thus, we only need to apply lemma \ref{35}.
\epr

Theorem \ref{55} is a solution to Bar-Hillel-Carnap paradox: tautologies are, in general, uninformative and contradictions are always without information. This is a consequence of the \emph{veridicality thesis} associated to the strong conception of semantic information: $\phi$ qualifies as semantic information only if $\phi$ is true (Cf. \cite{floridi2011}). The fact that some tautologies can be informative is a consequence of the veridicality thesis plus a \emph{complexity thesis}: $\phi$ qualifies as semantic information only if $\phi$ is complex. In the present paper, the complexity of some tautologies is associated to have new symbols. We can qualify even more the informativity of tautologies, but we will delegate that for future works. By now, we want to conclude with a strong solution to Hintikka's scandal of deduction.

\bte

For every collection of updates $\bar{\mathcal{D}}$ over an $S$-database $D$ and valid deduction $\phi_{0},\ldots,\phi_{n}$ over $S$ with support $\{\phi_{0},\ldots,\phi_{k}\}\{\phi_{n}\}$, if $\phi_{0},\ldots,\phi_{n}$ is over $S$, then $I_{\bar{\mathcal{D}}}((\phi_{0}\wedge\cdots\wedge\phi_{k})\to\phi_{n})=0$ but not necessarily $I_{\bar{\mathcal{D}}}(\{\phi_{0},\ldots,\phi_{k}\}\{\phi_{n}\})=0$. Moreover, if $\phi_{0},\ldots,\phi_{n}$ has symbols not in $S$, there are collections of updates $\bar{\mathcal{D}}$ over $\mathcal{D}$ for which $I_{\bar{\mathcal{D}}}((\phi_{0}\wedge\cdots\wedge\phi_{k})\to\phi_{n}) > 0$.

\ete
\bpr
Let $D$ be an $S$-database. Since $\{\phi_{0},\ldots,\phi_{k}\}\{\phi_{n}\}$ is the support of a valid deduction, $(\phi_{0}\wedge\cdots\wedge\phi_{k})\to\phi_{n}$ is a tautology. Then, by theorem \ref{53}, $I_{\bar{\mathcal{D}}}((\phi_{0}\wedge\cdots\wedge\phi_{k})\to\phi_{n}) = C_{\bar{\mathcal{D}}}((\phi_{0}\wedge\cdots\wedge\phi_{k})\to\phi_{n})$. Applying lemma \ref{35}, we have that either $I_{\bar{\mathcal{D}}}((\phi_{0}\wedge\cdots\wedge\phi_{k})\to\phi_{n})=0$ or $I_{\bar{\mathcal{D}}}((\phi_{0}\wedge\cdots\wedge\phi_{k})\to\phi_{n}) > 0$ according $(\phi_{0} \wedge \cdots \wedge \phi_{k})\to\phi_{n}$ has or not only symbols in the signature $S$. In any case, by lemma \ref{46}, $I_{\bar{\mathcal{D}}}(\{\phi_{1},\ldots,\phi_{k}\}\{\phi_{n}\})>0$ if for some $\mathcal{D}\in\bar{\mathcal{D}}$ we have that $\mathcal{D}(\phi_n)=1$ and at least one $\phi\in\{\phi_{1},\ldots,\phi_{k}\}$, $A(\phi)=0$, $\mathcal{D}(\phi)=1$ and $T\nvDash\phi$. Indeed, as $\mathcal{D}(\phi)=1$ and $T\nvDash\phi$, we have that $T$ is not a complete theory of $A$. Since $\phi\in\{\phi_{1},\ldots,\phi_{k}\}$, we know that $\{\phi_{0},\ldots,\phi_{k}\} \neq \oslash$. Finally, because $\mathcal{D}(\phi_n)=1$, $\phi_n$ is $D$-acceptable.
\epr

This theorem is our strong solution to the scandal of deduction. It is strong because it shows not only that some of the premisses of the scandal are false. It demonstrates that we can have informative deductions, notwithstanding the uninformativity of their associated conditionals. In other words, as Hintikka wanted to show, the conclusion of the scandal of deduction is false: there are informative deductions. This is, in turn, a consequence of the veridicality and complexity theses plus an \emph{inferential thesis}:  $\phi_{0},\ldots,\phi_{n}$ qualifies as semantic information only if $\phi_{0},\ldots,\phi_{n}$ is relevant.

\section{Conclusion}

We have proposed a way of measuring the semantic informativity of deductions by means of dynamic concepts of informational complexity and relevancy. In an schematic form, we can express our approach in the following way:

\[
\mbox{ Semantic informativity }=\mbox{ Complexity }\times\mbox{ Relevancy}.
\]

With respect to finite databases, we showed, first, how some tautologies can be informative but contradictions cannot. This is our solution to Bar-Hillel-Carnap paradox. It relies on two thesis: veridicality and complexity of semantic information. Adding one more thesis, inferentiality, we also obtained a solution to Hintikka' scandal of deductions: some deductions can be informative and others not. These three theses are part of the schematic form above.

As we have tried to make clear in the previous sections, this reformulation of Floridi's strong theory of semantic information poses some challenges to the veridicality thesis, but it does not invalidate it. In fact, the core is the veridicality thesis, which in some way characterize the left term of the schematic form above. What we have proposed in this paper is that we can measure the left term using the two terms of the right side of the equality above.

Given the work done here, we think there is a lot of future works to be done. In \cite{araujo2014}, we showed some consequences for artificial intelligence of the quantitative approach to semantic information developed in this paper, but there is other important directions to be explored. For applications, we can, for instance, to define databases in terms of typed constructive structures, and so to analyze the semantic information processing in some functional programming language, as it was made in \cite{primiero2013} with respect to data errors. Moreover, with respect to the foundations of the strong theory of semantic information, it would be desirable, for example, to implement levels of abstraction into databases. In this way, we will be able to evaluate Floridi's theses about the epistemic status of the levels of abstraction for semantic information. This is just to mention some important possible developments - the reader will find some other directions of research in \cite{araujo2014}.

\section*{Acknowledgements}

I would like to thank Viviane Beraldo for her support, to Luciano Floridi for his comments on my lecture at PT-AI2013 about the subject of this paper, and to Pedro Carrasqueira for his comments on a previous version of it. This research was funded by S\~ao Paulo Research Foundation (FAPESP) [2011/07781-2].

\bibliographystyle{spbasic}

\bibliography{bibliografia}

\end{document}